\documentclass[12pt,draftclsnofoot,onecolumn]{IEEEtran}

\usepackage{comment,url}
\usepackage[dvips]{graphicx}
\usepackage{caption,subfig}
\usepackage{float}
\graphicspath{{../fig/}}
\DeclareGraphicsExtensions{.eps}
\usepackage{setspace}
\usepackage{amsmath}
\usepackage{amssymb}
\usepackage{cite}
\usepackage{array}
\usepackage{mdwmath}
\usepackage{mdwtab}

\long\def\symbolfootnote[#1]#2{\begingroup%
\def\thefootnote{\fnsymbol{footnote}}\footnote[#1]{#2}\endgroup}

\begin{document}

\title{Random Distances Associated With\\Equilateral Triangles} 

\author{Yanyan Zhuang and Jianping Pan\\
University of Victoria, Victoria, BC, Canada}

\maketitle

\begin{abstract}
In this report, the explicit probability density functions of the random
Euclidean distances associated with equilateral triangles are given, when the
two endpoints of a link are randomly distributed in 1) the same triangle, 2)
two adjacent triangles sharing a side, 3) two parallel triangles sharing a 
vertex, and 4) two diagonal triangles sharing a vertex, respectively. The 
density function for 1) is based on the most recent work by 
B{\"a}sel~\cite{bäsel2012random}. 2)--4) are based on 1) and our previous 
work in~\cite{yanyan2011, yanyan2011h}. Simulation results
show the accuracy of the obtained closed-form distance distribution functions,
which are important in the theory of geometrical probability. The first two
statistical moments of the random distances and the polynomial fits of the
density functions are also given in this report for practical uses.
\end{abstract}


\section{The Problem}

\begin{figure}
  \centering
  \includegraphics[width=0.4\columnwidth]{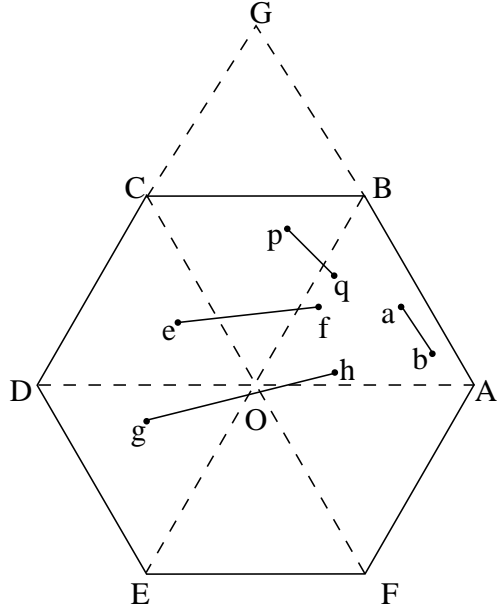}
  \caption{Random Points and Distances Associated with Equilateral Triangles.}
  \label{fig:triangle}
\end{figure}

Define a ``unit triangle'' as the equilateral triangle with side length $1$.
Picking two points uniformly at random from the interior of a unit triangle, or
between two adjacent unit triangles sharing a side or a vertex, the goal is to 
obtain the probabilistic density function (PDF) of the random distances between 
these two endpoints. 
There are four different cases $|ab|$, $|pq|$, $|ef|$ and $|gh|$, depending 
on the geometric locations of these two random
endpoints, as shown in Fig.~\ref{fig:triangle}. 
The next section gives the explicit PDFs for these cases.

\section{Distance Distributions Associated with Equilateral Triangles}\label{sec:result}

\subsection{$|ab|$: Distance Distribution within an Equilateral Triangle}
The author of~\cite{bäsel2012random} obtained the chord length distribution 
function for \textit{any} regular polygon. From this result, 
\cite{bäsel2012random} further derived the density function
for the distance between two uniformly and independently distributed
random points in the regular polygon. Although the methods used were  
elementary, this work can be considered as a major breakthrough in 
Geometrical Probability, which also helps us verify the distance distribution 
in a regular hexagon~\cite{yanyan2011h}.

Following are the notations used in~\cite{bäsel2012random}: $\mathcal P_{\rm
n,r}$ is the regular polygon with $n$ vertices and with a circumscribed circle of 
radius $r$; $l_k$ is the distance between vertices, given by
\begin{equation}
l_k=2r\sin \frac{k\pi}{n}, 
\end{equation}
where $k=0,1,...,K$ and $K=\lfloor\frac{n-2}{2}\rfloor$. $L$ denotes the perimeter 
and $A$ the area of $\mathcal P_{\rm n,r}$:
\begin{equation}
L=2nr\sin \frac{\pi}{n} ~~\mbox{ and }~~ A=\frac{1}{2}nr^2\sin \frac{2\pi}{n}.
\end{equation}

Denote the chord length distribution derived in~\cite{bäsel2012random} as
$f(s)$ for $\mathcal P_{\rm n,r}$, and the density function for the distance 
between two random points in $\mathcal P_{\rm n,r}$ as $g_{D_{\rm I}}(d)$, 
the relationship between these two 
functions is as follows according to~\cite{piefke1978beziehungen}:

\begin{equation}
g_{D_{\rm I}}(d)=\frac{2Ld}{A^2}\int_d^{l_{\rm K+1}}(s-d)f(s){\rm d}s.
\end{equation}

According to this relationship and the derived chord length distribution $f(s)$,
the density function of random distances in an equilateral triangle, $g_{D_{\rm I}}(d)$, 
is a special case in~\cite{bäsel2012random} when $n=3,r=\frac{1}{\sqrt{3}}$:
\begin{equation}\label{eq:triangle_within}
   g_{D_{\rm I}}(d)= 4d\left\{
     \begin{array}{lcr}

\left(2+\frac{4\pi}{3\sqrt{3}}\right)d^2-8d+\frac{2\pi}{\sqrt{3}} 
& \quad & 0\leq d\leq \frac{\sqrt{3}}{2}\\

\frac{2}{\sqrt{3}}\left(4d^2+6\right)\sin^{-1}\frac{\sqrt{3}}{2d}+\left(2-
\frac{8\pi}{3\sqrt{3}}\right)d^2+6\sqrt{4d^2-3} \\
~~~~-8d-\frac{4\pi}{\sqrt{3}}  & \quad & \frac{\sqrt{3}}{2}\leq d \leq 1\\

       0 & \quad & {\rm otherwise}
     \end{array}
   \right..
\end{equation}

The corresponding cumulative distribution function (CDF) is
\begin{equation}
   G_{D_{\rm I}}(d)= 2\left\{
     \begin{array}{lcr}
0 & \quad & d \leq 0 \\

\left(1+\frac{2\pi}{3\sqrt{3}}\right)d^4-\frac{16}{3}d^3+\frac{2\pi}{\sqrt{3}}
d^2 & \quad & 0\leq d\leq \frac{\sqrt{3}}{2} \\

\frac{4d^2}{\sqrt{3}}\left(d^2+3\right)\sin^{-1}\frac{\sqrt{3}}{2d}+\left(\frac{
26d^2}{3}+1\right)\sqrt{d^2-\frac{3}{4}}+\left(1-\frac{4\pi}{3\sqrt{3}}
\right)d^4 \\
~~~~-\frac{16}{3}d^3-\frac{4\pi}{\sqrt{3}}d^2  & \quad & \frac{\sqrt{3}}{2}\leq
d \leq 1\\

       1 & \quad & d \geq 1
     \end{array}
   \right..
\end{equation}

\subsection{$|pq|$: Distance Distribution between Two Adjacent Equilateral Triangles 
Sharing a Side}
Given the result above, and the result obtained by us in~\cite{yanyan2011} for the
distance distribution in a rhombus, further results of the random distances
between two equilateral triangles can be derived. In
Fig.~\ref{fig:triangle}, rhombus $OABC$ can be decomposed into two congruent,
adjacent equilateral triangles $\Delta OAB$ and $\Delta OBC$. Picking points
uniformly at random from the interior of this rhombus, then the points are 
equally likely to fall inside any one of these two triangles. Therefore, given the 
location of one endpoint of a random link, the second endpoint falls inside the same 
triangle as the first one with probability $\frac{1}{2}$ (such as $|ab|$), and with 
probability $\frac{1}{2}$ falls inside the adjacent triangle (such as $|pq|$). 

Suppose rhombus $OABC$ in Fig.~\ref{fig:triangle} has a side length of $1$, 
then the distribution of $|ab|$ is known as $g_{D_{\rm I}}(d)$ in 
(\ref{eq:triangle_within}) above. Denote the distribution of $|pq|$
as $g_{D_{\rm A}}(d)$.  The probability density function of the random 
distances between two uniformly distributed points that are both inside the 
same rhombus is $f_{D_{\rm I}}(d)$ (see (1) in~\cite{yanyan2011}). From the 
reasoning in the previous paragraph,
\begin{equation}
f_{D_{\rm I}}(d)=\frac{1}{2}g_{D_{\rm I}}(d)+\frac{1}{2}g_{D_{\rm A}}(d),
\end{equation}
and we have
\begin{equation}
g_{D_{\rm A}}(d)=2 f_{D_{\rm I}}(d)-g_{D_{\rm I}}(d).
\end{equation}

Therefore, the probability density function of the random distances
    between two uniformly distributed points, one in each of the two adjacent 
    unit triangles that are sharing a side, is
\begin{equation}\label{eq:triangle_btw}
  g_{D_{\rm A}}(d)=4d\left\{
    \begin{array}{lcr}

\frac{8}{3}d-\left(\frac{2}{3}+\frac{10\pi}{9\sqrt{3}}\right)d^2 
& \quad & 0\leq d\leq \frac{\sqrt{3}}{2}\\

-\frac{4}{\sqrt{3}}\left(1+\frac{4d^2}{3}\right)\sin^{-1}\frac{\sqrt{3}}{2d}
+\left(\frac{14\pi}{9\sqrt{3}}-\frac{2}{3}\right)d^2-\frac{8}{3}\sqrt{4d^2-3}
 \\
~~~~+\frac{8}{3}d+\frac{2\pi}{\sqrt{3}} & \quad & 
\frac{\sqrt{3}}{2}\leq d\leq 1\\

\frac{4}{\sqrt{3}}\left(1-\frac{d^2}{3}\right)\sin^{-1}\frac{\sqrt{3}}{2d}
+\left(\frac{2\pi}{9\sqrt{3}}-\frac{2}{3}\right)d^2+\sqrt{4d^2-3} \\
~~~~-\frac{2\pi}{3\sqrt{3}}-1 & \quad & 1\leq d\leq \sqrt{3} \\

      0 & \quad & {\rm otherwise}
    \end{array}
  \right..
\end{equation}

The corresponding CDF is

\begin{equation}\label{eq:triangle_btw_cdf}
   G_{D_{\rm A}}(d)= 2\left\{
     \begin{array}{lcr}
0 & \quad & d \leq 0 \\

\frac{16}{9}d^3-\left(\frac{1}{3}+\frac{5\pi}{9\sqrt{3}}\right)d^4 
& \quad & 0\leq d\leq \frac{\sqrt{3}}{2} \\

-\frac{4d^2}{\sqrt{3}}\left(1+\frac{2d^2}{3}\right)\sin^{-1}\frac{\sqrt{3}}{2d}
-4d^2\sqrt{d^2-\frac{3}{4}}+\left(\frac{7\pi}{9\sqrt{3}}-\frac{1}{3}\right)d^4
\\
~~~~+\frac{16}{9}d^3+\frac{2\pi}{\sqrt{3}}d^2  & \quad & \frac{\sqrt{3}}{2}\leq
d \leq 1\\

\frac{4d^2}{\sqrt{3}}\left(1-\frac{d^2}{6}\right)\sin^{-1}\frac{\sqrt{3}}{2d}
+\left(\frac{11d^2}{9}+\frac{5}{6}\right)\sqrt{d^2-\frac{3}{4}}+\left(\frac{\pi}
{9\sqrt{3}}-\frac{1}{3}\right)d^4 \\
~~~~-\left(\frac{2\pi}{3\sqrt{3}}+1\right)d^2-\frac{1}{4}  & \quad &
1\leq d \leq \sqrt{3}\\

       1 & \quad & d \geq \sqrt{3}
     \end{array}
   \right..
\end{equation}

Note that although unit triangles are assumed in
(\ref{eq:triangle_within})--(\ref{eq:triangle_btw_cdf}), the distance
distribution functions can be easily scaled by a nonzero scalar, for equilateral triangles
of arbitrary side length. For example, let the side length of such triangles be
$s>0$, then
\begin{equation}
G_{sD}(d)=P(sD\leq d)=P(D\leq \frac{d}{s})=G_D(\frac{d}{s}). \nonumber
\end{equation}
Therefore,
\begin{equation}\label{eq:scale}
 g_{sD}(d)=G'_D(\frac{d}{s})=\frac{1}{s}g_D(\frac{d}{s}).
\end{equation}

\subsection{$|ef|$: Distance Distribution between Two Parallel Equilateral Triangles 
Sharing a Vertex}

This case corresponds to the random distance $|ef|$ in Fig.~\ref{fig:triangle}. Here four 
unit triangles $\Delta OAB$, $\Delta OBC$, $\Delta OCD$ and $\Delta BCG$ together create a 
larger equilateral triangle $\Delta AGD$ with side length $2$. According to (\ref{eq:scale}), 
the density function of distance distribution inside triangle $\Delta AGD$ is 
$g_{2D_{\rm I}}(d)=\frac{1}{2}g_{D_{\rm I}}(\frac{d}{2})$, as $s=2$. On the other hand, if we 
look at the two random endpoints of a given link inside the large triangle, they will fall 
into one of the two following cases: i) one of the endpoints falls inside one of the three unit 
triangles on the boarder of the large triangle, such as $\Delta OAB$, $\Delta OCD$ or $\Delta BCG$, 
with probability $\frac{3}{4}$; ii) one of the endpoints falls inside the unit triangle 
$\Delta OBC$ in the middle, with probability $\frac{1}{4}$. 
Each of these two cases includes several more detailed sub-cases as follows:

\renewcommand{\labelenumi}{Case \roman{enumi})}
\begin{enumerate}
 \item Given the location of the first endpoint, the second endpoint will fall inside the same 
 triangle as the first one (such as $|ab|$) with probability $\frac{1}{4}$, fall inside 
 the adjacent triangle sharing a side (such as $|pq|$) with probability $\frac{1}{4}$, and 
 fall inside one of the parallel triangles sharing a vertex (such as $|ef|$) with probability 
 $\frac{1}{2}$. 

 \item When the location of the first endpoint is in $\Delta OBC$, the second endpoint 
 will fall inside the same triangle with probability $\frac{1}{4}$, and fall 
 inside one of the adjacent triangles sharing a side with probability $\frac{3}{4}$. 

\end{enumerate}

Denote the density function of random distance $|ef|$ as $g_{D_{\rm P}}(d)$, we have the 
following 
\begin{equation}
g_{2D_{\rm I}}(d)=\frac{3}{4}\left[\frac{1}{4}g_{D_{\rm I}}(d)+\frac{1}{4}g_{D_{\rm A}}(d)+
\frac{1}{2}g_{D_{\rm P}}(d)\right]+\frac{1}{4}\left[\frac{1}{4}g_{D_{\rm I}}(d)+
\frac{3}{4}g_{D_{\rm A}}(d)\right].
\end{equation}

Hence,
\begin{equation}
g_{D_{\rm P}}(d)=\frac{8}{3}g_{2D_{\rm I}}(d)-\frac{2}{3}g_{D_{\rm I}}(d)-g_{D_{\rm A}}(d)
=\frac{4}{3}g_{D_{\rm I}}(\frac{d}{2})-\frac{2}{3}g_{D_{\rm I}}(d)-g_{D_{\rm A}}(d).
\end{equation}

Therefore, the probability density function of the random distances between two uniformly
distributed points, one in each of the two parallel unit triangles that are sharing a 
vertex, is
\begin{equation}\label{eq:triangle_para}
  g_{D_{\rm P}}(d)=4d\left\{
    \begin{array}{lcr}

\left(\frac{4\pi}{9\sqrt{3}}-\frac{1}{3}\right)d^2 
& \quad & 0\leq d\leq \frac{\sqrt{3}}{2}\\

-\frac{4}{\sqrt{3}}\sin^{-1}\frac{\sqrt{3}}{2d}+\left(\frac{4\pi}{9\sqrt{3}}-
\frac{1}{3}\right)d^2-\frac{4}{3}\sqrt{4d^2-3}+\frac{2\pi}{\sqrt{3}} & \quad & 
\frac{\sqrt{3}}{2}\leq d\leq 1\\

\frac{4}{\sqrt{3}}\left(\frac{d^2}{3}-1\right)\sin^{-1}\frac{\sqrt{3}}{2d}
+d^2-\sqrt{4d^2-3}-\frac{8}{3}d+\frac{2\pi}{\sqrt{3}}+1 & \quad & 1\leq d\leq \sqrt{3} \\

\frac{4}{\sqrt{3}}\left(\frac{d^2}{3}+2\right)\sin^{-1}\frac{\sqrt{3}}{d}+\left(\frac{1}{3}
-\frac{4\pi}{9\sqrt{3}}\right)d^2+4\sqrt{d^2-3}\\
~~~~-\frac{8}{3}d-\frac{8\pi}{3\sqrt{3}} & \quad & \sqrt{3} \leq d \leq 2\\

      0 & \quad & {\rm otherwise}
    \end{array}
  \right..
\end{equation}

The corresponding CDF is

\begin{equation}\label{eq:triangle_para_cdf}
   G_{D_{\rm P}}(d)= 2\left\{
     \begin{array}{lcr}
0 & \quad & d \leq 0 \\

\left(\frac{2\pi}{9\sqrt{3}}-\frac{1}{6}\right)d^4 
& \quad & 0\leq d\leq \frac{\sqrt{3}}{2} \\

-\frac{4d^2}{\sqrt{3}}\sin^{-1}\frac{\sqrt{3}}{2d}+\left(\frac{2\pi}{9\sqrt{3}}
-\frac{1}{6}\right)d^4-\left(\frac{8d^2}{9}+\frac{1}{3}\right)\sqrt{4d^2-3}
+\frac{2\pi}{\sqrt{3}}d^2  & \quad & \frac{\sqrt{3}}{2}\leq d \leq 1\\

\frac{4d^2}{\sqrt{3}}\left(\frac{d^2}{6}-1\right)\sin^{-1}\frac{\sqrt{3}}{2d}
-\left(\frac{11d^2}{18}+\frac{5}{12}\right)\sqrt{4d^2-3}+\frac{d^4}{2}-\frac{16}{9}d^3 \\
~~~~+\left(\frac{2\pi}{\sqrt{3}}+1\right)d^2-\frac{1}{12}  & \quad &
1\leq d \leq \sqrt{3}\\

\frac{4d^2}{\sqrt{3}}\left(\frac{d^2}{6}+2\right)\sin^{-1}\frac{\sqrt{3}}{d}
+\left(\frac{26d^2}{9}+\frac{4}{3}\right)\sqrt{d^2-3}+\left(\frac{1}{6}
-\frac{2\pi}{9\sqrt{3}}\right)d^4 \\
~~~~-\frac{16}{9}d^3-\frac{8\pi}{3\sqrt{3}}d^2-\frac{5}{6} & \quad & 
\sqrt{3} \leq d \leq 2\\

       1 & \quad & d \geq 2
     \end{array}
   \right..
\end{equation}

\subsection{$|gh|$: Distance Distribution between Two Diagonal Equilateral Triangles 
Sharing a Vertex}

This case corresponds to the random distance $|gh|$ in Fig.~\ref{fig:triangle}. Here
a regular hexagon is divided into six unit triangles. Looking at the two random 
endpoints of a given link inside the hexagon, the first endpoint can fall inside any one 
of the six triangles, and the second endpoint will i) fall inside the same triangle as the 
first one (such as $|ab|$) with probability $\frac{1}{6}$; ii) fall inside the adjacent 
triangle sharing a side (such as $|pq|$) with probability $\frac{1}{3}$; iii) fall inside 
the parallel triangle sharing a vertex (such as $|ef|$) with probability $\frac{1}{3}$; 
iv) fall inside the diagonal triangle sharing a vertex (such as $|gh|$) with 
probability $\frac{1}{6}$. 

The density function of the random distances within a regular hexagon has been derived 
in~\cite{yanyan2011h}, and we denote it as $f_H(d)$. Also denote the density function 
of random distance $|gh|$ as $g_{D_{\rm D}}(d)$, we have
\begin{equation}
f_H(d)=\frac{1}{6}g_{D_{\rm I}}(d)+\frac{1}{3}g_{D_{\rm A}}(d)+\frac{1}{3}g_{D_{\rm P}}(d)
+\frac{1}{6}g_{D_{\rm D}}(d),
\end{equation}
or, 
\begin{equation}
g_{D_{\rm D}}(d)=6f_H(d)-\left[g_{D_{\rm I}}(d)+2g_{D_{\rm A}}(d)+2g_{D_{\rm P}}(d) 
\right]. 
\end{equation}

Therefore, the probability density function of the random distances between two 
uniformly distributed points, one in each of the two diagonal unit triangles that are 
sharing a vertex, is

\begin{equation}\label{eq:triangle_diag}
  g_{D_{\rm D}}(d)=4d\left\{
    \begin{array}{lcr}

\left(\frac{2}{3}-\frac{2\pi}{9\sqrt{3}}\right)d^2 
& \quad & 0\leq d\leq \frac{\sqrt{3}}{2}\\

\frac{4}{\sqrt{3}}\left(\frac{2d^2}{3}+1\right)\sin^{-1}\frac{\sqrt{3}}{2d}+
\left(\frac{2}{3}-\frac{14\pi}{9\sqrt{3}}\right)d^2+2\sqrt{4d^2-3}-\frac{2\pi}{\sqrt{3}} 
& \quad & \frac{\sqrt{3}}{2}\leq d\leq 1\\

-\frac{4}{\sqrt{3}}\left(\frac{2d^2}{3}+1\right)\sin^{-1}\frac{\sqrt{3}}{2d}
+\left(\frac{2\pi}{9\sqrt{3}}-\frac{2}{3}\right)d^2-2\sqrt{4d^2-3}\\
~~~~+\frac{16}{3}d+\frac{2\pi}{3\sqrt{3}} & \quad & 1\leq d\leq \sqrt{3} \\

-\frac{4d^2}{3\sqrt{3}}\sin^{-1}\frac{\sqrt{3}}{d}+\left(\frac{4\pi}{9\sqrt{3}}
-\frac{4}{3}\right)d^2-\frac{4}{3}\sqrt{d^2-3}+\frac{16}{3}d-4 & \quad & 
\sqrt{3} \leq d \leq 2\\

      0 & \quad & {\rm otherwise}
    \end{array}
  \right..
\end{equation}

The corresponding CDF is
\begin{equation}\label{eq:triangle_diag_cdf}
   G_{D_{\rm D}}(d)= 2\left\{
     \begin{array}{lcr}
0 & \quad & d \leq 0 \\

\left(\frac{1}{3}-\frac{\pi}{9\sqrt{3}}\right)d^4 & \quad & 0\leq d\leq 
\frac{\sqrt{3}}{2} \\

\frac{4d^2}{\sqrt{3}}\left(\frac{d^2}{3}+1\right)\sin^{-1}\frac{\sqrt{3}}{2d}
+\left(\frac{13d^2}{9}+\frac{1}{6}\right)\sqrt{4d^2-3}+\left(\frac{1}{3}-
\frac{7\pi}{9\sqrt{3}}\right)d^4\\
~~~~-\frac{2\pi}{\sqrt{3}}d^2  & \quad & \frac{\sqrt{3}}{2}\leq d \leq 1\\

-\frac{4d^2}{\sqrt{3}}\left(\frac{d^2}{3}+1\right)\sin^{-1}\frac{\sqrt{3}}{2d}
-\left(\frac{13d^2}{9}+\frac{1}{6}\right)\sqrt{4d^2-3}+\left(\frac{\pi}{9\sqrt{3}}
-\frac{1}{3}\right)d^4 \\
~~~~+\frac{32}{9}d^3+\frac{2\pi}{3\sqrt{3}}d^2+\frac{1}{3}  
& \quad & 1\leq d \leq \sqrt{3}\\

-\frac{2d^4}{3\sqrt{3}}\sin^{-1}\frac{\sqrt{3}}{d}+\left(\frac{4}{3}-
\frac{10d^2}{9}\right)\sqrt{d^2-3}+\left(\frac{2\pi}{9\sqrt{3}}
-\frac{2}{3}\right)d^4+\frac{32}{9}d^3 \\
~~~~-4d^2+\frac{11}{6} & \quad & \sqrt{3} \leq d \leq 2\\

       1 & \quad & d \geq 2
     \end{array}
   \right..
\end{equation}

\section{Verification by Simulation}

\begin{figure}
  \centering
  \includegraphics[width=0.6\columnwidth]{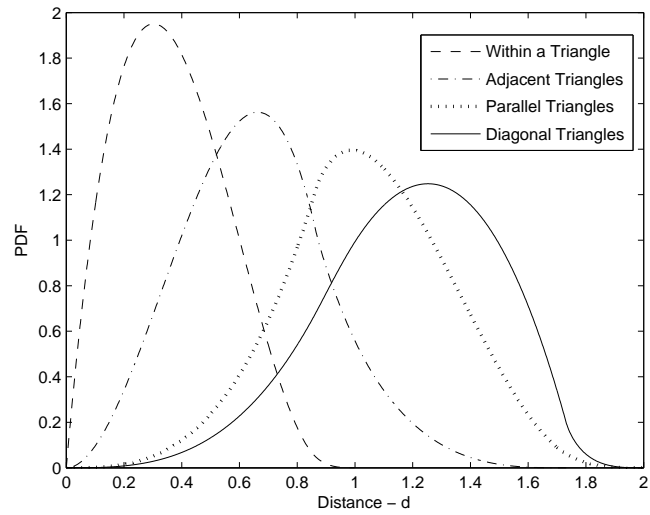}
  \caption{Distributions of Random Distances Associated with Equilateral Triangles.}
  \label{fig:triangle_pdf}
\end{figure}

\begin{figure}
  \centering
  \includegraphics[width=0.6\columnwidth]{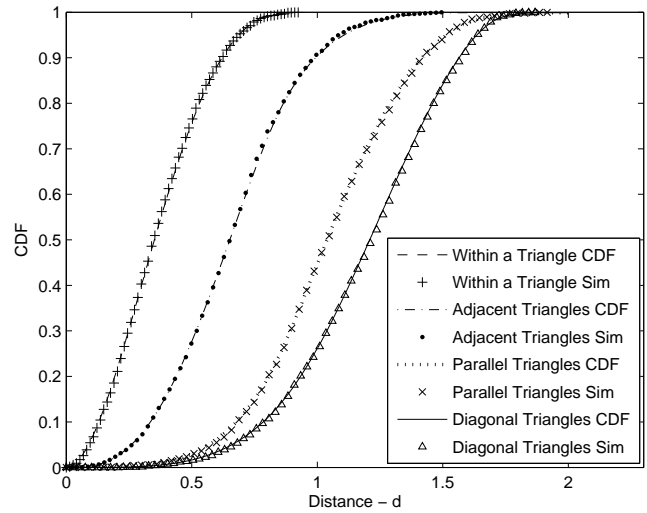}
  \caption{Distribution and Simulation Results for Random Distances Associated
with Equilateral Triangles.}
  \label{fig:triangle_cdf}
\end{figure}

Figure~\ref{fig:triangle_pdf} plots the probability density functions, as given
in (\ref{eq:triangle_within}), (\ref{eq:triangle_btw}), (\ref{eq:triangle_para}) 
and (\ref{eq:triangle_diag}) of the four random distance cases shown in 
Fig.~\ref{fig:triangle}. Figure~\ref{fig:triangle_cdf} 
shows a comparison between the cumulative distribution functions (CDFs) of the 
random distances, and the simulation results by generating $1,000$ pairs of 
random points with the corresponding geometric locations. 
Figure~\ref{fig:triangle_cdf} demonstrates that our distance distribution 
functions are very accurate when compared with the simulation results.

\section{Practical Results}

\subsection{Statistical Moments of Random Distances}

The distance distribution functions given in Section~\ref{sec:result} can
conveniently lead to all the statistical moments of the random distances
associated with equilateral triangles. Given $g_{D_{\rm I}}(d)$ in
(\ref{eq:triangle_within}), for example, the first moment (mean) of $d$, i.e.,
the average distance within a unit triangle, is 
\begin{eqnarray}
 M_{D_{\rm I}}^{(1)}=\int_0^1xg_{D_{\rm I}}(x)dx=\frac{1}{5}+\frac{3}{20}\ln
(3)\approx 0.3647918433, \nonumber
\end{eqnarray}
and the second raw moment is
\begin{eqnarray}
 M_{D_{\rm I}}^{(2)}=\int_0^1x^2g_{D_{\rm I}}(x)dx=\frac{1}{6},
\nonumber
\end{eqnarray}
from which the variance (the second central moment) can be derived as
\begin{eqnarray}
 Var_{D_{\rm I}}=M_{D_{\rm I}}^{(2)}-\left[M_{D_{\rm I}}^{(1)}\right]^2\approx
0.0335935777. \nonumber
\end{eqnarray}

When the side length of the unit triangle is scaled by $s$, the corresponding first 
two statistical moments given above then become
\begin{equation}
M_{D_{\rm I}}^{(1)}=0.3647918433s,~~\mbox{}~~M_{D_{\rm I}}^{(2)}=\frac{s}{6}
~~\mbox{ and }~~Var_{D_{\rm I}}=0.0335935777s^2.
\end{equation}

\begin{table}
  \caption{Moments and Variance---Numerical vs Simulation Results}
  \centering
  \begin{tabular}{|c||c|c|c|c|}
    \hline
    Endpoint Geometry & PDF/Sim & $M_{D}^{(1)}$ & $M_{D}^{(2)}$ & $Var_{D}$ \\
\hline \hline
    Within a & $g_{D_{\rm I}}(d)$ & $0.3647918433s$ & $0.1666666667s$ &
$0.0335935777s^2$ \\ 
    \cline{2-5}
    Single Triangle & Sim & $0.3636606517s$ & $0.1654245278s$ &
$0.0331164521s^2$ \\ \hline
    Between Two & $g_{D_{\rm A}}(d)$ & $0.6599648287s$ & $0.5s$ &
$0.0644464249s^2$\\ 
    \cline{2-5}
    Adjacent Triangles & Sim & $0.6597703174s$ & $0.4991102154s$ &
$0.0637684635s^2$ \\ \hline
Between Two & $g_{D_{\rm P}}(d)$ & $1.0423971067s$ & $1.1666666667s$ &
$0.0800749386s^2$\\ 
    \cline{2-5}
    Parallel Triangles & Sim & $1.0423259678s$ & $1.1655336144s$ &
$0.0791204824s^2$ \\ \hline
Between Two & $g_{D_{\rm D}}(d)$ & $1.1880379828s$ & $1.5s$ &
$0.0885657513s^2$\\ 
    \cline{2-5}
    Diagonal Triangles & Sim & $1.1896683303s$ & $1.5048850531s$ &
$0.0896042750s^2$ \\ \hline
  \end{tabular}
  \label{tab:moment}
\end{table}

Table~\ref{tab:moment} lists the first two moments and the variance of the
random distances in all four cases given in Section~\ref{sec:result}. It also
gives the corresponding simulation results for verification purposes. 

\subsection{Polynomial Fits of Random Distances}

\begin{table}
  \caption{Coefficients of the Polynomial Fit and the Norm of Residuals (NR)}
  \centering
  \begin{tabular}{|c||c|c|}
    \hline
    PDF & Polynomial Coefficients & NR \\ \hline \hline 
     \vspace*{-0.12cm} & \vspace*{-0.12cm} & \vspace*{-0.12cm} \\
     & $10^{10}\times
\left[0.006410~-0.062752~~0.283869~-0.787308~~1.497828\right.$ & \\
     $g_{D_{\rm I}}(d)$ &
$-2.071998~~2.155623~-1.720734~~1.065808~-0.514668~~0.193640$ & $0.002646$ \\
     & $-0.056449~~0.012613~-0.002124~~0.000263~ -0.000023~~0.000001$ & \\
     & $\left.0~~0~0~~0\right]$ & \\\hline
     \vspace*{-0.12cm} & \vspace*{-0.12cm} & \vspace*{-0.12cm} \\
& $10^8\times
\left[-0.000192~~0.003514~-0.029768~~0.154523~-0.549699\right.$ & \\  
    $g_{D_{\rm A}}(d)$ &
$1.420142~-2.754830~~4.091879~-4.703977~~4.202832~-2.914966$ &  \\ 
     & $1.559704~-0.636497~~0.194663~-0.043503~~0.006952~-0.000722$ & $0.086864$
\\ 
     & $\left.0.000047~-0.000002~~0~0\right]$ & \\\hline
\vspace*{-0.12cm} & \vspace*{-0.12cm} & \vspace*{-0.12cm} \\
& $10^7\times
\left[-0.000069~~0.001559~-0.016145~~0.101938~-0.439383\right.$ & \\  
    $g_{D_{\rm P}}(d)$ &
$1.370859~-3.202401~~5.714218~-7.873972~~8.415340~-6.697807$ &  \\ 
     & $4.441839~-2.155120~~0.781941~-0.206851~~0.038480~-0.004776$ & $0.105657$
\\ 
     & $\left.0.000365~-0.000015~~0~0\right]$ & \\\hline
     \vspace*{-0.12cm} & \vspace*{-0.12cm} & \vspace*{-0.12cm} \\
& $10^7\times
\left[0.000023~-0.000530~~0.005579~-0.035274~~0.150448\right.$ & \\  
    $g_{D_{\rm D}}(d)$ &
$-0.460164~~1.045986~-1.805145~~2.394131~-2.453508~~1.942627$ &  \\ 
     & $-1.182240~~0.547275~-0.189550~~0.047945~-0.008552~~0.001022$ & $0.075017$
\\ 
     & $\left.0.000076~-0.000003~~0~0\right]$ & \\\hline
  \end{tabular}
  \label{tab:poly}
\end{table}

\begin{figure}
\centering
  \subfloat[Within a Single
Triangle]{\includegraphics[width=0.5\columnwidth]{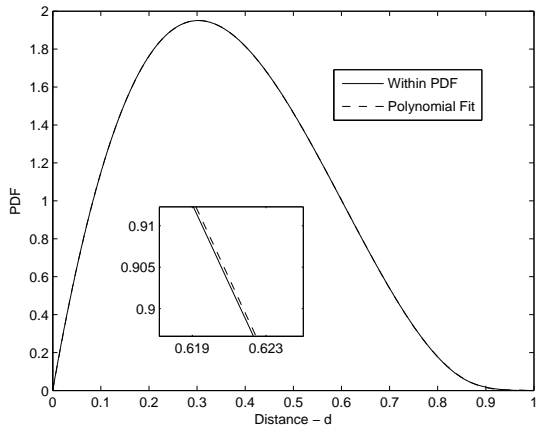}}
  \subfloat[Between two Adjacent Triangles Sharing a 
  Side]{\includegraphics[width=0.5\columnwidth]{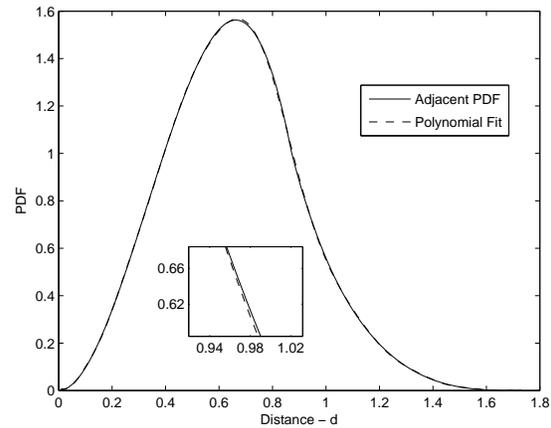}}\\
\subfloat[Between two Parallel Triangles Sharing a 
  Vertex]{\includegraphics[width=0.5\columnwidth]{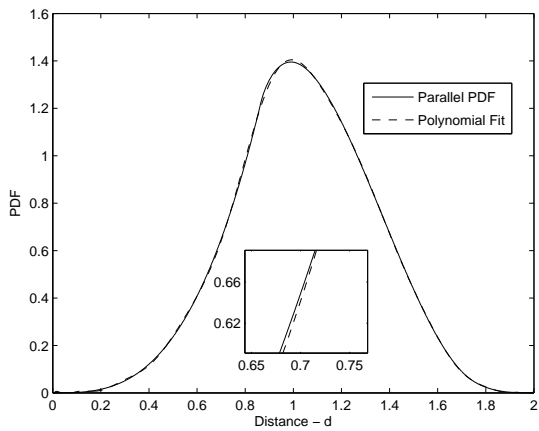}}
  \subfloat[Between two Diagonal Triangles Sharing a 
  Vertex]{\includegraphics[width=0.5\columnwidth]{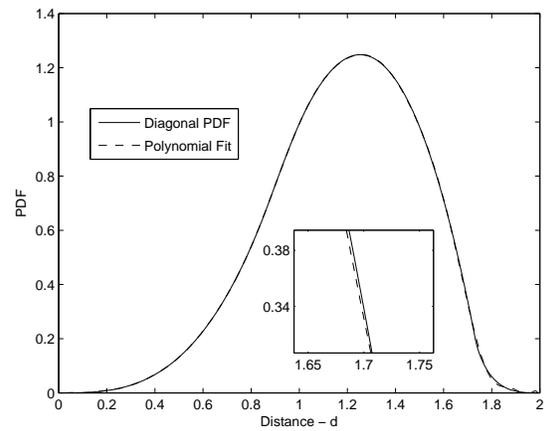}}
  \caption{Polynomial Fit of the Distance Distribution Functions Associated with
Equilateral Triangles.}
  \label{fig:triangle_poly}
\end{figure}

Table~\ref{tab:poly} lists the coefficients of the degree-$20$ polynomial fits
of the original PDFs given in Section~\ref{sec:result}, from $d^{20}$ to
$d^{0}$, and the corresponding norm of residuals. 
Figure~\ref{fig:triangle_poly}(a)--(d) plot the polynomials listed in 
Table~\ref{tab:poly} with the original PDFs. From the
figure, it can be seen that all the polynomials match closely with the original
PDFs. These high-order polynomials facilitate further manipulations of the
distance distribution functions, with a high accuracy. 

\section{Conclusions}
\label{sec:conclude}

In this report, we gave the closed-form probability density functions of the
random distances associated with equilateral triangles. The correctness of the
obtained results has been validated by simulation. The first two statistical
moments and the polynomial fits of the density functions are also given for
practical uses. 

\section*{Acknowledgment}
The authors would like to thank 
Dr. Aaron Gulliver for initially posing the problem associated with
hexagons, which lead to our previous work~\cite{yanyan2011, 
yanyan2011h} and this report.

\end{document}